\RequirePackage{ifpdf}
\ifpdf 
\documentclass[pdftex]{sigma}
\else
\documentclass{sigma}
\fi

\newcommand\ZZ{\mathbb{Z}}

\newcommand\FSA{{\cal A}}

\newcommand\al\alpha
\newcommand\be\beta
\newcommand\ga\gamma
\newcommand\la\lambda

\newcommand\half{\frac12}
\newcommand\thalf{\tfrac12}
\newcommand\const{{\rm const}\,}
\newcommand{\qhyp}[5]{\,\mbox{}_{#1}\phi_{#2}\left(
  \genfrac{}{}{0pt}{}{#3}{#4};#5\right)}
\newcommand\LHS{left-hand side}

\newcommand\goH{\mathfrak{H}}
\newcommand\DAHA{\tilde\goH}
\newcommand\Dsym{D_{\rm sym}}
\newcommand\Zsym{Z_{\rm sym}}
\newcommand\Asym{\FSA_{\rm sym}}
\newcommand\AW{\widetilde{AW}(3)}
\newcommand\AWQ{\widetilde{AW}(3,Q_0)}
\newcommand\Q{\widetilde Q}

\numberwithin{equation}{section}

\newtheorem{Theorem}{Theorem}[section]
\newtheorem{Lemma}[Theorem]{Lemma}
\newtheorem{Proposition}[Theorem]{Proposition}
\newtheorem{Corollary}[Theorem]{Corollary}
\newtheorem{Remark}[Theorem]{Remark}
\newtheorem{Definition}[Theorem]{Definition}

\begin{document}

\allowdisplaybreaks

\renewcommand{\PaperNumber}{063}

\FirstPageHeading

\renewcommand{\thefootnote}{$\star$}

\ShortArticleName{Zhedanov's Algebra $AW(3)$}

\ArticleName{The Relationship between Zhedanov's Algebra $\boldsymbol{AW(3)}\!$
and the Double Af\/f\/ine Hecke Algebra\\ in the Rank One Case\footnote{This paper
is a contribution to the Vadim Kuznetsov Memorial Issue
`Integrable Systems and Related Topics'. The full collection is
available at
\href{http://www.emis.de/journals/SIGMA/kuznetsov.html}{http://www.emis.de/journals/SIGMA/kuznetsov.html}}}

\Author{Tom H. KOORNWINDER}

\AuthorNameForHeading{T.H. Koornwinder}

\Address{Korteweg-de Vries Institute, University of Amsterdam,\\
Plantage Muidergracht 24, 1018 TV Amsterdam, The Netherlands}
\Email{\href{mailto:thk@science.uva.nl}{thk@science.uva.nl}}
\URLaddress{\url{http://www.science.uva.nl/~thk/}}

\ArticleDates{Received December 22, 2006, in f\/inal form April
23, 2007; Published online April 27, 2007}

\Abstract{Zhedanov's algebra $AW(3)$
is considered with explicit structure constants
such that, in the basic representation,
the f\/irst generator becomes the second order
$q$-dif\/ference operator for the Askey--Wilson polynomials.
It is proved that this
representation is faithful for a certain quotient of $AW(3)$
such that the Casimir operator is equal to
a special constant. Some explicit aspects of the double
af\/f\/ine Hecke algebra (DAHA) related to symmetric and non-symmetric Askey--Wilson
polynomials are presented and proved without requiring knowledge of
general DAHA theory. Finally a central extension of this quotient of
$AW(3)$ is introduced
which can be embedded in the DAHA by means of the faithful basic
representations of both algebras.}


\Classification{33D80; 33D45}

\begin{flushright}
{\em Dedicated to the memory of Vadim Kuznetsov}
\end{flushright}

\noindent
{\sl Briefly after I had moved from CWI, Amsterdam to a professorship at
the University of Amsterdam in 1992, Vadim Kuznetsov contacted me about
the possibility to come to Amsterdam as a postdoc. We successfully applied
for a grant. He arrived with his wife Olga and his son Simon in Amsterdam
for a two-years stay during 1993--1995. I vividly remember picking them up
at the airport and going in the taxi with all their stuff to their first
apartment in Amsterdam, at the edge of the red light quarter. These were
two interesting years, where we learnt a lot from each other. We wrote one
joint paper, but Vadim wrote many further papers alone or with other coauthors
during this period. We should have written more together, but our temperaments
were too different for that. Vadim was always speeding up, while I wanted
to ponder and to look for further extensions and relations with other work.

After his Amsterdam years Vadim had a marvelous career which led to
prestigious UK grants, tenure in Leeds, and a lot of organizing of
conferences and proceedings. We met several times afterwards. I visited for
instance Leeds for one week, and Vadim was an invited speaker
at the conference in Amsterdam in 2003 on the occasion of my sixtieth birthday.}

\section{Introduction}
Zhedanov \cite{1} introduced in 1991 an algebra $AW(3)$
with three generators $K_0$, $K_1$, $K_2$ and
three relations in the form of $q$-commutators, which describes deeper
symmetries of the Askey--Wilson polynomials. In fact,
for suitable choices of the structure constants of the algebra, the
Askey--Wilson polynomial $p_n(x)$ is the kernel of an
intertwining operator between a representation of $AW(3)$ by
$q$-dif\/ference operators on the space of polynomials in $x$ and
a representation by tridiagonal operators on the space of
of inf\/inite sequences $(c_n)_{n=1,2,\ldots}$. In the f\/irst representation
$K_1$ is multiplication by $x$ and $K_0$ is the second order
$q$-dif\/ference operator for which the Askey--Wilson polynomials are
eigenfunctions with explicit eigenvalues $\la_n$.
In the second representation $K_0$ is the diagonal operator with diagonal
elements $\la_n$ and $K_1$ is the tridiagonal operator corresponding to
the three-term recurrence relation for the Askey--Wilson polynomials.
The formula for $p_n(x)$ expressing the intertwining property with respect
to $K_2$ is the so-called {\em $q$-structure relation} for the
Askey--Wilson polynomials (see \cite{6}) and the relation for $AW(3)$
involving the $q$-commutator of $K_1$ and $K_2$ is the so-called
{\em $q$-string equation} (see \cite{7}).
Terwilliger \& Vidunas \cite{16} showed that every Leonard pair satisf\/ies
the $AW(3)$ relations for a suitable choice of the structure constants.

In 1992, one year after Zhedanov's paper \cite{1},
Cherednik \cite{8} introduced double af\/f\/ine Hecke algebras associated
with root systems (DAHA's). This was the f\/irst of an important series
of papers by the same author, where a representation of the DAHA was
given in terms of $q$-dif\/ference-ref\/lection operators
($q$-analogues of Dunkl operators), joint eigenfunctions
of such operators were identif\/ied as non-symmetric Macdonald polynomials,
and Macdonald's conjectures for ordinary (symmetric) Macdonald polynomials
associated with root systems could be proved. For a nice exposition
of this theory see Macdonald's recent book \cite{5}.
In particular, the DAHA approach to Macdonald--Koornwinder polynomials,
due to several authors
(see Sahi~\cite{9,10}, Stokman \cite{14} and references given there)
is also presented in~\cite{5}.
The last chapter of~\cite{5} discusses the rank one specialization
of these general results. For the DAHA of type $A_1$ (one parameter)
this yields non-symmetric $q$-ultraspherical polynomials.
For the DAHA of type $(C_1^\vee,C_1)$ (four parameters)
the non-symmetric Askey--Wilson
polynomials are obtained. These were earlier treated by Sahi~\cite{10} and by
Noumi \& Stokman \cite{11}.
See also Sahi's recent paper~\cite{12}.

Comparison of Zhedanov's $AW(3)$ with the DAHA of type
of type $(C_1^\vee,C_1)$, denoted by $\DAHA$,
suggests some relationship. Both algebras are
presented by generators and relations, the f\/irst has a representation
by $q$-dif\/ference operators on the space of symmetric Laurent polynomials
in~$z$ and the second has a representation by
$q$-dif\/ference-ref\/lection operators on the space of general Laurent polynomials
in $z$. Since this representation of the DAHA is called the
{\em basic representation} of $\DAHA$, I will call the just mentioned
representation of $AW(3)$ also the {\em basic representation}.
In the basic representation of $AW(3)$ the operator $K_0$ is
equal to some operator $D$ occurring in the basic representation of $\DAHA$
and involving ref\/lections, provided
$D$ is restricted in its action to symmetric Laurent polynomials.
This suggests that the basic representation of $AW(3)$ may remain
valid if
we represent $K_0$ by $D$, so that it involves ref\/lection terms.
It will turn out in this paper that this conjecture is correct in the
$A_1$ case, i.e., when the Askey--Wilson parameters are restricted to
the continuous $q$-ultraspherical case. In the general case the conjecture
is true for a rather harmless central extension of $AW(3)$ involving
a generator $T_1$, which will be identif\/ied with the familiar $T_1$
in $\DAHA$ which has in the basic representation of $\DAHA$ the
symmetric Laurent polynomials as one of its two eigenspaces.

This paper does not suppose any knowledge about the general theory of
double af\/f\/ine Hecke algebras and about Macdonald and related polynomials in
higher rank.
The contents of the paper
are as follows. Section 2 presents $AW(3)$ and
its relationship with Askey--Wilson polynomials.
We add to $AW(3)$ one more relation expressing that the
Casimir operator $Q$ is equal to a~special constant $Q_0$
(of course precisely the constant occurring for $Q$ in the basic
representation),
and we denote the resulting quotient algebra
by $AW(3,Q_0)$. Then it is shown
that the basic representation of $AW(3,Q_0)$ is faithful.
Section~3 discusses~$\DAHA$ (the DAHA of type $(C_1^\vee,C_1)$), its basic
representation, and the basis vectors
for the 2-dimensional eigenspaces of the operator $D$
in terms of Askey--Wilson polynomials.
Section 4 gives an explicit expression for the non-symmetric
Askey--Wilson polynomials which is in somewhat dif\/ferent terms than the explicit
expression in \cite[\S~6.6]{5}.
Two presentations of $\DAHA$ by generators and relations of PBW-type
are given in Section 5. The faithfulness of the basic representation is
proved (a result which of course is also a special case of the known
result in the case of general rank, see Sahi~\cite{9}).
The main result of the present paper, the embedding
of a central extension of $AW(3,Q_0)$ in
$\DAHA$, is stated and proved in Section~6.

For the computations in this paper I made heavy use of computer algebra
performed in
{\em Mathematica}${}^{\mbox{\footnotesize\textregistered}}$.
For reductions of
expressions in non-commuting variables subject to relations I used
the package {\em NCAlgebra} \cite{13} within
{\em Mathematica}${}^{\mbox{\footnotesize\textregistered}}$.
{\em Mathematica} notebooks containing these computations will be available
for downloading in
\url{http://www.science.uva.nl/~thk/art/}.

\bigskip

\noindent
{\bf Conventions}\\
Throughout assume that $q$ and $a$, $b$, $c$, $d$ are complex constants such that
\begin{gather}
q\ne0,\qquad
q^m\ne1\  (m=1,2,\ldots),\qquad
a,b,c,d\ne0,\qquad
abcd\ne q^{-m}\  (m=0,1,2,\ldots).\!\!
\label{21}
\end{gather}
Let $e_1$, $e_2$, $e_3$, $e_4$ be the elementary symmetric polynomials in $a$, $b$, $c$, $d$:
\begin{gather}
e_1:=a+b+c+d,\qquad
e_2:=ab+ac+bc+ad+bd+cd,\nonumber\\
e_3:=abc+abd+acd+bcd,\qquad
e_4:=abcd.
\label{58}
\end{gather}
For \mbox{($q$-)Pochhammer} symbols and
\mbox{($q$-)hypergeometric series}
use the notation of \cite{2}. In particular,
\begin{gather*}
(a;q)_k:=\prod_{j=0}^{k-1}(1-aq^j),\qquad
(a_1,\ldots,a_r;q)_k:=(a_1;q)_k\cdots(a_r;q)_k,
\\
\qhyp r{r-1}{q^{-n},a_2,\ldots,a_r}{b_1,\ldots,b_{r-1}}{q,z}:=
\sum_{k=0}^n\frac{(q^{-n},a_2,\ldots,a_r;q)_k}
{(b_1,\ldots,b_{r-1},q;q)_k}\,z^k.
\end{gather*}
For Laurent polynomials $f$
in $z$ the $z$-dependence will be written as $f[z]$.
Symmetric Laurent polynomials
$f[z]=\sum\limits_{k=-n}^n c_k z^k$ (where $c_k=c_{-k}$) are related to ordinary
polynomials $f(x)$ in $x=\thalf(z+z^{-1})$ by
$f(\thalf(z+z^{-1}))=f[z]$.

\section[Zhedanov's algebra $AW(3)$]{Zhedanov's algebra $\boldsymbol{AW(3)}$}
\label{5}

Zhedanov \cite{1} introduced an algebra $AW(3)$ with three generators
$K_0$, $K_1$, $K_2$ and with three relations
\begin{gather*}
[K_0,K_1]_q=K_2,\\
[K_1,K_2]_q=B\,K_1+ C_0\,K_0+D_0,\\
[K_2,K_0]_q=B\,K_0+C_1\,K_1+D_1,
\end{gather*}
where
\[
[X,Y]_q:=q^\half XY-q^{-\half} YX
\]
is the $q$-commutator and
where the structure constants
$B$, $C_0$, $C_1$, $D_0$, $D_1$ are f\/ixed complex constants.
He also gave a {\em Casimir operator}
\begin{gather*}
Q:=\big(q^{-\half}-q^{\frac32}\big)K_0K_1K_2+qK_2^2+B(K_0K_1+K_1K_0)+qC_0K_0^2
+q^{-1}C_1K_1^2\\
\phantom{Q:=}{}+(1+q)D_0K_0+(1+q^{-1})D_1K_1,
\end{gather*}
which commutes with the generators.

Clearly, $AW(3)$
can equivalently be described as an algebra with two generators
$K_0$, $K_1$ and with two relations
\begin{gather}
\label{1}
(q+q^{-1})K_1K_0K_1-K_1^2K_0-K_0K_1^2=B\,K_1+ C_0\,K_0+D_0,\\
\label{2}
(q+q^{-1})K_0K_1K_0-K_0^2K_1-K_1K_0^2=B\,K_0+C_1\,K_1+D_1.
\end{gather}
Then the Casimir operator $Q$ can be written as
\begin{gather}
Q=(K_1K_0)^2\!-(q^2+1+q^{-2})K_0(K_1K_0)K_1+(q+q^{-1})K_0^2K_1^2\!
+(q+q^{-1})(C_0K_0^2+C_1K_1^2)\nonumber\\
\phantom{Q=}{}+B\bigl((q+1+q^{-1})K_0K_1+K_1K_0\bigr)
+(q+1+q^{-1})(D_0K_0+D_1K_1).\label{80}
\end{gather}

Let the structure constants be expressed
in terms of $a$, $b$, $c$, $d$ by means of $e_1$, $e_2$, $e_3$, $e_4$
(see~\eqref{58}) as follows:
\begin{gather}
B:=(1-q^{-1})^2(e_3+qe_1),\nonumber\\
C_0:=(q-q^{-1})^2,\nonumber\\
C_1:=q^{-1}(q-q^{-1})^2 e_4,\label{42}\\
D_0:=-q^{-3}(1-q)^2(1+q)(e_4+qe_2+q^2),\nonumber\\
D_1:=-q^{-3}(1-q)^2(1+q)(e_1e_4+qe_3).\nonumber
\end{gather}
Then there is a representation (the {\em basic representation}) of the algebra
$AW(3)$ with structure constants \eqref{42}
on the space $\Asym$
of symmetric Laurent polynomials $f[z]=f[z^{-1}]$
as follows:
\begin{gather}
(K_0f)[z]=(\Dsym f)[z],\qquad
(K_1f)[z]=((Z+Z^{-1})f)[z]:=(z+z^{-1})f[z],
\label{4}
\end{gather}
where $\Dsym$, given by
\begin{gather}
(\Dsym f)[z]:=\frac{(1-az)(1-bz)(1-cz)(1-dz)}{(1-z^2)(1-qz^2)}\,
\bigl(f[qz]-f[z]\bigr)\nonumber\\
\phantom{(\Dsym f)[z]:=}{}+\frac{(a-z)(b-z)(c-z)(d-z)}{(1-z^2)(q-z^2)}\,
\bigl(f[q^{-1}z]-f[z]\bigr)
+(1+q^{-1}abcd)f[z],
\label{3}
\end{gather}
is the second order operator having the
{\em Askey--Wilson polynomials}
(see \cite{3}, \cite[\S~7.5]{2}, \cite[\S~3.1]{4})
as eigenfunctions.
It can indeed be verif\/ied that the operators $K_0$, $K_1$ given by
\eqref{4} satisfy relations \eqref{1}, \eqref{2} with structure
constants \eqref{42}, and that the Casimir operator $Q$ becomes the
following constant in this representation:
\begin{gather}
(Qf)(z)=Q_0\,f(z),
\label{82}
\end{gather}
where
\begin{gather}
Q_0:=q^{-4}(1-q)^2\Bigl(q^4(e_4-e_2)+q^3(e_1^2-e_1e_3-2e_2)\nonumber\\
\phantom{Q_0:=}{}-q^2(e_2e_4+2e_4+e_2)
+q(e_3^2-2e_2e_4-e_1e_3)+e_4(1-e_2)\Bigr).
\label{81}
\end{gather} 

Let $AW(3,Q_0)$ be the algebra generated by $K_0$, $K_1$ with relations
\eqref{1}, \eqref{2} and
\begin{gather}
Q=Q_0,
\label{83}
\end{gather}
assuming the structure constants \eqref{42}.
Then the basic representation of $AW(3)$ is also a representation of
$AW(3,Q_0)$.

The Askey--Wilson polynomials are given by
\begin{gather}
p_n\bigl(\thalf(z+z^{-1});a,b,c,d\mid q\bigr)
:=\frac{(ab,ac,ad;q)_n}{a^n}\,
\qhyp43{q^{-n},q^{n-1}abcd,az,az^{-1}}{ab,ac,ad}{q,q}.
\label{11}
\end{gather}
These polynomials are symmetric in $a$, $b$, $c$, $d$ (although this cannot be read
of\/f from \eqref{11}).
We will work with the renormalized version which is {\em monic} as a
Laurent polynomial in $z$ (i.e., the coef\/f\/icient of $z^n$ equals 1):
\begin{gather}
P_n[z]=P_n[z;a,b,c,d\mid q]:=
\frac1{(abcdq^{n-1};q)_n}\,
p_n\bigl(\thalf(z+z^{-1});a,b,c,d\mid q\bigr)
\nonumber\\
\phantom{P_n[z]}{}=a^{-n}\sum_{k=0}^n
\frac{(q^{-n};q)_k\,(az,az^{-1};q)_k\,(abq^k,acq^k,adq^k;q)_{n-k}\,q^k}
{(q;q)_k\,(abcdq^{n+k-1};q)_{n-k}}\,.
\label{13}
\end{gather}
Note that the monic
Askey--Wilson polynomials $P_n[z]$ are well-def\/ined for all $n$ under
condition~\eqref{21}.

The eigenvalue equation involving $\Dsym$ is
\begin{gather}
\Dsym P_n=\la_nP_n,\qquad
\la_n:=q^{-n}+abcd q^{n-1}.
\label{12}
\end{gather}
Under condition \eqref{21} all eigenvalues in \eqref{12} are distinct.

The three-term recurrence relation for the monic Askey--Wilson polynomials
(see \cite[(3.1.5)]{4}) is as follows:
\begin{gather}
(z+z^{-1})P_n[z]=P_{n+1}[z]+\be_n P_n[z]+\ga_n P_{n-1}[z]\qquad(n\ge1),\nonumber\\
(z+z^{-1})P_0[z]=P_1[z]+\be_0 P_0[z],
\label{59}\\
\be_n:=q^{n-1}\,\frac{(1-q^n-q^{n+1})e_3+qe_1+q^{2n-1}e_3e_4
-q^{n-1}(1+q-q^{n+1})e_1e_4}{(1-q^{2n-2}e_4)(1-q^{2n}e_4)},
\label{60}\\
\ga_n:=(1-q^{n-1}ab)(1-q^{n-1}ac)(1-q^{n-1}ad)(1-q^{n-1}bc)
(1-q^{n-1}bd)(1-q^{n-1}cd)
\nonumber\\
\phantom{\ga_n:=}{}\times\frac{(1-q^n)(1-q^{n-2}e_4)}
{(1-q^{2n-3}e_4)(1-q^{2n-2}e_4)^2(1-q^{2n-1}e_4)}.
\label{61}
\end{gather}
From this we see that $P_n[z]$ remains well-def\/ined if the condition
$a,b,c,d\ne0$ in \eqref{21}
is omitted.
It also follows from \eqref{12} and \eqref{59}--\eqref{61} that the
representation \eqref{4} of $AW(3)$ is not necessarily irreducible, but that
it has $1\in\Asym$ as a cyclic element. The representation
will become irreducible if we moreover require that none of
$ab$, $ac$, $ad$, $bc$, $bd$, $cd$ equals $q^{-m}$ for some $m=0,1,2,\ldots$.

We now show that $AW(3,Q_0)$ has the elements
\begin{gather}
K_0^n(K_1K_0)^lK_1^m\qquad(m,n=0,1,2,\ldots, \ \ l=0,1)
\label{54}
\end{gather}
as a basis and that the representation \eqref{4} of $AW(3,Q_0)$ is faithful.
\begin{Lemma}
\label{55}
Each element of $AW(3,Q_0)$ can be written as a linear combination of ele\-ments \eqref{54}.
\end{Lemma}
\begin{proof}
$AW(3,Q_0)$ is spanned by elements $K_\al=K_{\al_1}\cdots K_{\al_k}$, where
$\al=(\al_1,\ldots,\al_k)$,
$\al_i=0$ or~1. Let $\rho(\al)$ the number of pairs $(i,j)$ such that $i<j$,
$\al_i=1$, $\al_j=0$.
$K_\al$ has the form~\eqref{54} if\/f $\rho(\al)=0$ or 1. We will show
that each $K_\al$ with $\rho(\al)>1$ can be written as a~linear combination
of elements $K_\be$ with $\rho(\be)<\rho(\al)$.
Indeed,
if $\rho(\al)>1$ then $K_\al$ must have a~substring $K_1K_1K_0$ or $K_1K_0K_0$
or $K_1K_0K_1K_0$. By substitution of relations
\eqref{1}, \eqref{2} or \eqref{83}
(with \eqref{80}), respectively, we see that each such string is a linear
combination of elements~$K_\be$ with~$\rho(\be)<\rho(\al)$.
\end{proof}
\begin{Theorem}
The elements \eqref{54} form a basis of $AW(3,Q_0)$ and the representation
\eqref{4} of $AW(3,Q_0)$ is faithful.
\label{76}
\end{Theorem}
\begin{proof}
Because of Lemma \ref{55} it is suf\/f\/icient to show that the operators
\begin{gather}
(\Dsym)^n\,(Z+Z^{-1})^m\qquad(m,n=0,1,2,\ldots),\nonumber\\
(\Dsym)^{n-1}\,(Z+Z^{-1})\,\Dsym\,(Z+Z^{-1})^{m-1}\qquad(m,n=1,2,\ldots)
\label{56}
\end{gather}
acting on $\Asym$ are linearly independent.
By \eqref{12} and \eqref{59} we have for all $j$:
\begin{gather}
(\Dsym)^n\,(Z+Z^{-1})^m\,P_j[z]
=\la_{j+m}^n P_{j+m}[z]+\cdots,\nonumber\\
(\Dsym)^{n-1}\,(Z+Z^{-1})\,\Dsym\,(Z+Z^{-1})^{m-1}\,P_j[z]
=\la_{j+m}^{n-1}\la_{j+m-1} P_{j+m}[z]+\cdots,
\label{84}
\end{gather}
where the right-hand sides give expansions in terms of $P_k[z]$ with
$k$ running from $j+m$ downwards.

Suppose that the operators \eqref{56} are not linearly independent. Then
\begin{gather}
\sum_{k=0}^m\sum_l a_{k,l} (\Dsym)^l\,(Z+Z^{-1})^k\nonumber\\
\qquad{}+
\sum_{k=1}^m\sum_l b_{k,l} (\Dsym)^{l-1}\,(Z+Z^{-1})\,\Dsym\,(Z+Z^{-1})^{k-1}
=0\label{85}
\end{gather}
for certain coef\/f\/icients $a_{k,l}$, $b_{k,l}$ such that for some $l$
$a_{m,l}\ne0$ or $b_{m,l}\ne0$. Then it follows from~\eqref{84}
that for all $j$, when we let the \LHS\ of \eqref{85} act on $P_j[z]$,
the coef\/f\/icient of $P_{j+m}[z]$ yields:
\begin{equation}
\sum_l (a_{m,l}\la_{j+m}^l+b_{m,l}\la_{j+m}^{l-1}\la_{j+m-1})=0.\label{86}
\end{equation}
By \eqref{12} we have, writing $x=q^{j+m}$ and $u= q^{-1}abcd$,
\[
\la_{j+m}=x^{-1}+ ux,\qquad
\la_{j+m-1}=qx^{-1}+ q^{-1}ux.
\]
We can consider the identity \eqref{86} as an
identity for Laurent polynomials in $x$.
Since the left-hand side vanishes
for infinitely many values of $x$,
it must be identically zero.
Let $n$ be the maximal $l$ for which $a_{m,l}\ne0$ or $b_{m,l}\ne0$.
Then, in particular, the coef\/f\/icients of $x^{-n}$ and $x^n$
in the
left-hand side of \eqref{86} must be zero.
This gives explicitly:
\begin{gather*}
a_{m,n}+q b_{m,n}=0,\qquad
u^n a_{m,n}+q^{-1}u^n b_{m,n}=0.
\end{gather*}
This implies $a_{m,n}=b_{m,n}=0$, contradicting our assumption.
\end{proof}
\begin{Remark}\rm
Note that we have 6 structure constants
$B$, $C_0$, $C_1$, $D_0$, $D_1$, $Q_0$ depending on~4~parameters $a$, $b$, $c$, $d$.
However, 2 degrees of freedom in the structure coef\/f\/icients are caused
by scale transformations. Indeed, the scale transformations
$K_0\to c_0K_0$ and $K_1\to c_1K_1$ induce the following transformations
on the structure coef\/f\/icients:
\begin{gather*}
B\to c_0c_1B,\quad
C_0\to c_1^2C_0,\quad
C_1\to c_0^2C_1,\quad
D_0\to c_0c_1^2D_0,\quad
D_1\to c_0^2c_1D_1,\quad
Q_0\to c_0^2c_1^2 Q_0.\!
\end{gather*}
But these scale transformations also af\/fect the basic representation.
This becomes $K_0=c_0\Dsym$, $K_1=c_1(Z+Z^{-1})$.
\end{Remark}

\section[The double affine Hecke algebra of type $(C_1^\vee,C_1)$]{The double af\/f\/ine Hecke algebra of type $\boldsymbol{(C_1^\vee,C_1)}$}

Recall condition \eqref{21}.
The double af\/f\/ine Hecke algebra of type $(C_1^\vee,C_1)$, denoted by
$\DAHA$
(see \cite[\S~6.4]{5}),
is generated by $Z$, $Z^{-1}$, $T_1$, $T_0$ with relations
$ZZ^{-1}=1=Z^{-1}Z$ and
\begin{gather}
(T_1+ab)(T_1+1)=0,
\label{6}\\
(T_0+q^{-1}cd)(T_0+1)=0,
\label{8}\\
(T_1Z+a)(T_1Z+b)=0,
\label{7}\\
(qT_0Z^{-1}+c)(qT_0Z^{-1}+d)=0.
\label{9}
\end{gather}
Here I have used the notation of \cite{12}, which is slightly
dif\/ferent from
the notation in \cite[\S~6.4]{5}. Conditions on $q$, $a$, $b$, $c$, $d$ in
\cite{5} are more strict than in \eqref{21}. This will give no problem,
as can be seen by checking all results hereafter from scratch.

From \eqref{6} and \eqref{8} and the non-vanishing of
$a$, $b$, $c$, $d$ we see that $T_1$ and $T_0$ are invertible:
\begin{gather}
T_1^{-1}=-a^{-1}b^{-1}T_1-(1+a^{-1}b^{-1}),
\label{38}\\
T_0^{-1}=-qc^{-1}d^{-1}T_0-(1+qc^{-1}d^{-1}).
\label{39}
\end{gather}
Put
\begin{gather}
Y:=T_1T_0,
\label{36}\\
D:=Y+q^{-1}abcdY^{-1}=T_1T_0+q^{-1}abcdT_0^{-1}T_1^{-1},
\label{37}\\
\Zsym:=Z+Z^{-1}.
\label{48}
\end{gather}
By \eqref{6} and \eqref{8} $D$ commutes with $T_1$ and $T_0$.
By \eqref{6} and \eqref{7} $\Zsym$ commutes with $T_1$.

The algebra $\DAHA$ has a faithful representation, the so-called
{\em basic representation}, on the space~$\FSA$
of Laurent polynomials $f[z]$ as follows:
\begin{gather}
(Zf)[z]:=z\,f[z],
\label{14}\\
(T_1f)[z]:=\frac{(a+b)z-(1+ab)}{1-z^2}\,f[z]+
\frac{(1-az)(1-bz)}{1-z^2}\,f[z^{-1}],
\label{15}\\
(T_0f)[z]:=\frac{q^{-1}z((cd+q)z-(c+d)q)}{q-z^2}\,f[z]
-\frac{(c-z)(d-z)}{q-z^2}\,f[qz^{-1}].
\label{16}
\end{gather}
The representation property is from \cite[\S~6.4]{5} or by straightforward
computation. The faithfulness is from \cite[(4.7.4)]{5} or by an independent
proof later in this paper.

Now we can compute:
\begin{gather}
(Yf)[z]=
\frac{z \bigl(1+ab-(a+b)z\bigr)
\bigl((c+d)q-(cd+q)z\bigr)}{q(1-z^2)(q-z^2)}\,f[z]\nonumber\\
\phantom{(Yf)[z]=}{}+\frac{(1-az)(1-bz)(1-cz)(1-dz)}{(1-z^2)(1-q z^2)}f[qz]\nonumber\\
\phantom{(Yf)[z]=}{}+\frac{(1-a z)(1-b z) \bigl((c+d)qz-(cd+q)\bigr)}
{q(1-z^2)(1-q z^2)}\,f[z^{-1}]\nonumber\\
\phantom{(Yf)[z]=}{}+\frac{(c-z)(d-z)\bigl(1+ab-(a+b)z\bigr)}{(1-z^2)(q-z^2)}\,f[qz^{-1}],
\label{34}
\\
(Df)[z]
=
\frac{(1-q)z(1-az) (1-bz)\,\bigl((q+1)(cd+q)z-q(c+d)(1+z^2)\bigr)}
{q(1-z^2)(q-z^2)(1-q z^2)}\,f[z^{-1}]
\nonumber\\
\phantom{(Df)[z]=}{}+\frac{(1-q)z(c-z)(d-z)\,\bigl((a+b) (q+z^2)-(a b+1)(q+1)z\bigr)}
{(1-z^2)(q-z^2)(q^2-z^2)}\,f[qz^{-1}]
\nonumber\\
\phantom{(Df)[z]=}{}+\Bigl((a+b)(cd+q)(q+z^2)+q(ab+1)(c+d)(1+z^2)
\nonumber\\
\phantom{(Df)[z]=}{}-\bigl((q+1)(cd+q)(ab+1)+2q(a+b)(c+d)\bigr)z\Bigr)
\,\frac z{q(1-z^2)(q-z^2)}\,f(z)
\label{10}\\
\phantom{(Df)[z]=}{}+\frac{(c-z)(d-z)(aq-z)(bq-z)\!}{(q-z^2)(q^2-z^2)} f[q^{-1}z]\!
+\frac{(1-az) (1-bz) (1-cz)(1-dz)\!}{(1-z^2)(1-qz^2)} f[qz]
.\!\nonumber
\end{gather}
If we compare \eqref{10} and \eqref{3} then we see that
\[
(Df)[z]=(\Dsym f)[z]\qquad\mbox{if}\qquad f[z]=f[z^{-1}].
\]
In particular, if we apply $D$ to the Askey--Wilson polynomial
$P_n[z]$ given by \eqref{13} then we obtain from \eqref{12} that
\begin{gather}
DP_n=\la_nP_n.
\label{17}
\end{gather}

By \eqref{6} and \eqref{8} the operators $T_1$ and $T_0$, acting on
$\FSA$ as given by \eqref{15}, \eqref{16} have two eigenvalues.
We can characterize the eigenspaces.
\begin{Proposition}
\label{19}
$T_1$ given by \eqref{15} has eigenvalues $-ab$ and $-1$.
$T_1f=-ab\,f$ iff $f$ is symmetric.
If $a$, $b$ are distinct from $a^{-1}$, $b^{-1}$ then
$T_1f=-f$ iff $f[z]=z^{-1}(1-az)(1-bz)g[z]$ for some symmetric Laurent
polynomial $g$.
\end{Proposition}

\begin{proof}
We compute
\begin{gather*}
(T_1f)[z]+ab\,f[z]=
\frac{(1-az)(1-bz)}{1-z^2}\,(f[z^{-1}]-f[z]),
\end{gather*}
which settles the f\/irst assertion. We also compute
\begin{gather*}
(T_1f)[z]+f[z]=
\frac{(1-az)(1-bz)}{1-z^2}\,f[z^{-1}]-
\frac{(a-z)(b-z)}{1-z^2}\,f[z].
\end{gather*}
This equals zero if $f[z]=z^{-1}(1-az)(1-bz)g[z]$  with $g$ symmetric.
On the other hand, if $(T_1f)[z]+f[z]=0$ and $a$, $b$ are distinct from
$a^{-1}$, $b^{-1}$ then
\[
(1-az)(1-bz)f[z^{-1}]=(a-z)(b-z)f[z]
\]
and hence $f[z]=z^{-1}(1-az)(1-bz)g[z]$ for some Laurent polynomial $g$ and
we obtain $g[z]=g[z^{-1}]$.
\end{proof}

\begin{Proposition}
\label{20}
$T_0$ given by \eqref{16} has eigenvalues $-q^{-1}cd$ and $-1$.
$T_0f=-q^{-1}cd\,f$ iff $f[z]=f[qz^{-1}]$.
If $c$, $d$ are distinct from $qc^{-1}$, $qd^{-1}$ then
$T_0f=-f$ iff $f[z]=z^{-1}(c-z)(d-z)g[z]$ for some Laurent
polynomial $g$ satisfying $g[z]=g[qz^{-1}]$.
\end{Proposition}

\begin{proof}
We compute
\begin{gather*}
(T_0f)[z]+q^{-1}cd\,f[z]=
\frac{(c-z)(d-z)}{q-z^2}\,(f[z]-f[qz^{-1}]),
\end{gather*}
which settles the f\/irst assertion. We also compute
\begin{gather*}
(T_0f)[z]+f[z]=
\frac{(q-cz)(q-dz)}{q(q-z^2)}\,f[z]-
\frac{q(c-z)(d-z)}{q(q-z^2)}\,f[qz^{-1}].
\end{gather*}
Then the second assertion is proved by similar arguments  as in the proof
of Proposition \ref{19}.\
\end{proof}

We now look for further explicit solutions of the eigenvalue equation
\begin{gather}
Df=\la_nf.
\label{18}
\end{gather}
Clearly, the solution $P_n$ (see \eqref{17}) also satisf\/ies
$T_1P_n=-ab\,P_n$. In order to f\/ind further solutions of \eqref{18} we make
an Ansatz for $f$ as suggested by Propositions \ref{19} and \ref{20},
namely $f[z]=z^{-1}(1-az)(1-bz)g[z]$ or
$f[z]=g[q^{-\half} z]$ or $f[z]=z^{-1}(c-z)(d-z)g[q^{-\half}z]$, in each
case with $g$ symmetric. Then it turns out that \eqref{18} takes the
form of the Askey--Wilson second order $q$-dif\/ference equation, but
with parameters and sometimes also the degree changed. We thus obtain
as
further solutions $f$ of \eqref{18} for $n\ge1$:
\begin{gather}
Q_n[z]:=a^{-1}b^{-1}z^{-1}(1-az)(1-bz)\,
P_{n-1}[z;qa, qb,c,d\mid q],
\label{31}\\
P_n^\dagger[z]:=q^{\half n}\,
P_n\big[q^{-\half}z;q^\half a,q^\half b,q^{-\half}c,q^{-\half}d\mid q\big],
\label{32}\\
Q_n^\dagger[z]:=q^{\half(n-1)}z^{-1}(c-z)(d-z)\,
P_{n-1}\big[q^{-\half}z;q^\half a,q^\half b,q^\half c,q^\half d\mid q\big].
\label{33}
\end{gather}
So we have for $n\ge1$
four dif\/ferent eigenfunctions of $D$ at
eigenvalue $q^{-n}+abcdq^{n-1}$ which are also eigenfunction of $T_1$ or $T_0$:
\begin{gather}
T_1P_n=-ab\,P_n,\qquad
T_1Q_n=-Q_n,\qquad
T_0P_n^\dagger=-q^{-1}cd\,P_n^\dagger,\qquad
T_0Q_n^\dagger=-Q_n^\dagger.
\label{51}
\end{gather}
They all are Laurent polynomials of degree $n$ with highest term
$z^n$ and lowest term
$\const z^{-n}$:
\begin{alignat}{3}
&P_n[z]=z^n+\cdots+z^{-n},\qquad&&
Q_n[z]=z^n+\cdots+a^{-1}b^{-1}z^{-n},&\nonumber\\
&P_n^\dagger[z]=z^n+\cdots+q^nz^{-n},\qquad&&
Q_n^\dagger[z]=z^n+\cdots+q^{n-1}cdz^{-n}.&
\label{29}
\end{alignat}
Since the eigenvalues $\la_n$ are distinct for dif\/ferent $n$,
it follows that $D$ has a 1-dimensional eigenspace $\FSA_0$ at
eigenvalue $\la_0$, consisting of the constant Laurent polynomials,
and that it has a
2-dimensional eigenspace $\FSA_n$ at eigenvalue $\la_n$ if $n\ge1$,
which has $P_n$ and $P_n^\dagger$ as basis vectors, but which also has
any other two out of $P_n$, $Q_n$, $P_n^\dagger$, $Q_n^\dagger$ as basis vectors,
provided these two functions have the coef\/f\/icients of $z^{-n}$ distinct.
Generically we can use any two out of these four as basis vectors.
The basis consisting of $P_n$ and $P_n^\dagger$ occurs in \cite[\S~6.6]{5}.
In the following sections we will work f\/irst with the basis consisting of
$P_n$ and $Q_n^\dagger$, but afterwards it will be more convenient to
use $P_n$ and $Q_n$.

\section[Non-symmetric Askey-Wilson polynomials]{Non-symmetric Askey--Wilson polynomials}

Since $T_1$ and $T_0$ commute with $D$, the eigenspaces of $D$ in $\FSA$
are invariant under $Y=T_1T_0$. We can f\/ind explicitly the eigenvectors of
$Y$ within these eigenspaces $\FSA_n$.
\begin{Theorem}
\label{79}
The non-symmetric Askey--Wilson polynomials
\begin{gather}
E_{-n}[z]:=\frac1{1-q^{n-1}cd}\,(P_n[z]-Q_n^\dagger[z])\qquad(n=1,2,\ldots),
\label{23}\\
E_n[z]:=\frac{q^n(1-q^{n-1}abcd)}{1-q^{2n-1}abcd}\,P_n[z]+
\frac{1-q^n}{1-q^{2n-1}abcd}\,Q_n^\dagger[z]\qquad(n=1,2,\ldots),
\label{24}\\
E_0[z]:=1
\label{35}
\end{gather}
span the one-dimensional eigenspaces of $Y$ within $\FSA_n$ with
the following eigenvalues:
\begin{gather}
YE_{-n}=q^{-n}\,E_{-n}\qquad(n=1,2,\ldots),
\label{25}\\
YE_n=q^{n-1}abcd\,E_n\qquad(n=0,1,2,\ldots).
\label{26}
\end{gather}
The coefficients of highest and lowest terms in $E_{-n}$ and $E_n$ are:
\begin{gather}
E_{-n}[z]=z^{-n}+\cdots+\const z^{n-1}\qquad(n=1,2,\ldots),
\label{27}
\\
E_n[z]=
z^n+\cdots+\left(1-\frac{(1-q^n)(1-q^{n-1}cd)}{1-q^{2n-1}abcd}\right)z^{-n}
\qquad(n=1,2,\ldots).
\label{28}
\end{gather}
\end{Theorem}
\begin{proof}
Clearly, by their def\/inition, $E_{-n}$ and $E_n$ are in $\FSA_n$,
while \eqref{27}, \eqref{28} follow from \eqref{29}.
Equation \eqref{26} for $n=0$ follows from \eqref{36} and Propositions
\ref{19} and \ref{20}.
For the proof of \eqref{25}, \eqref{26} we use a
$q$-dif\/ference equation for Askey--Wilson polynomials
(see \cite[(7.7.7)]{2}, \cite[(3.1.8)]{4}):
\begin{gather}
\frac{P_n[q^{-\half}z;a,b,c,d\mid q]-P_n[q^\half z;a,b,c,d\mid q]}
{(q^{-\half n}-q^{\half n})(z-z^{-1})}=
P_{n-1}\big[z;q^\half a,q^\half b,q^\half c,q^\half d\mid q\big].
\label{30}
\end{gather}
The expression $(YE_{-n})[z]-q^{-n}E_{-n}[z]$ ($n=1,2,\ldots$)
only involves terms
$P_n[w;a,b,c,d\mid q]$ for $w=z,qz,q^{-1}z$ and
terms $P_{n-1}[w;q^\half a,q^\half b,q^\half c,q^\half d\mid q]$
for $w=q^{-\half} z, q^\half z$, as can be seen from \eqref{23}, \eqref{33}
and \eqref{34}. Now twice substitute in this expression \eqref{30} with
$z$ replaced by $q^{-\half} z$ and $q^\half z$, respectively. Then we arrive
at an expression only involving terms $P_n[w;a,b,c,d\mid q]$
for $w=z,qz,q^{-1}z$. By \eqref{3} it can be recognized as
$((\Dsym P_n)[z]-(q^{-n}+abcd q^{n-1})P_n[z])/(1-q^n)$,
which equals zero
by \eqref{12}. This settles \eqref{25}.
The reduction of the expression
$(YE_n)[z]-q^{n-1}abcd E_n[z]$ ($n=1,2,\ldots$)
can be done in a completely similar way. Here we arrive at
the expression
$((\Dsym P_n)[z]-(q^{-n}+abcd q^{n-1})P_n[z])
/(1-q^{1-2n}(abcd)^{-1})$, which equals zero.
\end{proof}
\begin{Remark}\rm
By condition \eqref{21} all eigenvalues of $Y$ on $\FSA$
(see \eqref{25}, \eqref{26}) are distinct. So for all $n\in\ZZ$
$E_n[z]$ is the unique Laurent polynomial of degree $|n|$ which satisf\/ies
\eqref{25} or \eqref{26} and has coef\/f\/icient of $z^n$ equal to 1.
Moreover, for $n\ge1$, $E_{-n}$ is the unique element of $\FSA_n$
of the form \eqref{27}, and $E_n$ is the unique element of $\FSA_n$
of the form \eqref{28}
\end{Remark}
\begin{Remark}\rm
The occurrence of the $q$-dif\/ference equation \eqref{30} in the proof of
Theorem \ref{79} and the occurrence of Askey--Wilson polynomials with shifted
parameters as eigenfunctions of $D$ (see \eqref{31}--\eqref{33})
is probably much related to the one-variable case of the
$q$-dif\/ference equations in Rains \cite[Corollary 2.4]{15}.
\end{Remark}

From \eqref{29}, \eqref{27} and \eqref{28} we obtain
\begin{gather}
E_{-n}=\frac{ab}{ab-1}\,(P_n-Q_n)\qquad(n=1,2,\ldots),
\label{49}\\
E_n=\frac{(1-q^n ab)(1-q^{n-1}abcd)}{(1-ab)(1-q^{2n-1}abcd)}\,P_n-
\frac{ab(1-q^n)(1-q^{n-1}cd)}{(1-ab)(1-q^{2n-1}abcd)}\,
Q_n\qquad(n=1,2,\ldots).\!\!\!
\label{50}
\end{gather}
Next, \eqref{49}, \eqref{50} and \eqref{51} yield
\begin{gather}
T_1 E_{-n}=-\frac{1+ab-abcdq^{n-1}-abq^n}{1-abcdq^{2n-1}}\,E_{-n}-ab\,E_n
\qquad(n=1,2,\ldots),
\label{52}\\
T_1 E_n=\frac{(1-q^n)(1-abq^n)(1-cdq^{n-1})(1-abcdq^{n-1})}
{(1-abcdq^{2n-1})^2}\,E_{-n}
\nonumber\\
\phantom{T_1 E_n=}{}
-\frac{abq^{n-1}(cd+q-cdq^n-abcdq^n)}{1-abcdq^{2n-1}}\,E_n\qquad
(n=1,2,\ldots).
\label{53}
\end{gather}

\section[A PBW-type theorem for $\tilde{\mathfrak{H}}$]{A PBW-type theorem for $\boldsymbol{\DAHA}$}

In this section I will give two other sets of relations for $\DAHA$, both
equivalent to \eqref{6}--\eqref{9} and both of PBW-type form.
For the second set of relations we will see that the spanning set of elements
of $\DAHA$, as implied by these relations, is indeed a basis. This is done
by showing that this set of elements is linearly independent in the basic
representation,  which also shows that this representation is faithful.
The faithfulness of the basic representation was f\/irst shown,
in the more general $n$ variable setting, by Sahi~\cite{9}.

\begin{Proposition}
$\DAHA$ can equivalently be described as the algebra generated by
$T_1$, $T_0$, $Z$, $Z^{-1}$ with relations
$ZZ^{-1}=1=Z^{-1}Z$ and
\begin{gather}
T_1^2=-(ab+1)T_1-ab,
\label{64}\\
T_0^2=-(q^{-1}cd+1)T_0-q^{-1}cd,
\label{65}\\
T_1Z =Z^{-1}T_1+(ab+1) Z^{-1}-(a+b),
\label{66}\\
T_1Z^{-1}=ZT_1-(ab+1)Z^{-1}+(a+b),
\label{67}\\
T_0Z=qZ^{-1}T_0-(q^{-1}cd+1)Z+(c+d),
\label{68}\\
T_0Z^{-1}=qZT_0+q^{-1}(q^{-1}cd+1)Z-q^{-1}(c+d).
\label{69}
\end{gather}
$\DAHA$ is spanned by the elements $Z^mT_0^iY^nT_1^j$, where
$m\in\ZZ$, $n=0,1,2,\ldots$, $i,j=0,1$.
\end{Proposition}

\begin{proof}
\eqref{64}, \eqref{66} are equivalent to \eqref{6}, \eqref{7}, and
\eqref{65}, \eqref{68} are equivalent to \eqref{8}, \eqref{9}.
Furthermore, \eqref{66} is equivalent to \eqref{67}, and
\eqref{68} is equivalent to \eqref{69}.
Hence relations \mbox{\eqref{64}--\eqref{69}} are equivalent to
relations \eqref{6}--\eqref{9}.

For the second statement note that \eqref{64}--\eqref{69} imply that
each word in $\DAHA$ can be written as a linear combination of
words $Z^mT_0^i(T_1T_0)^nT_1^j$, where
$m\in\ZZ$, $n=0,1,2,\ldots$, $i,j=0,1$. Then substitute $Y=T_1T_0$.
\end{proof}
\begin{Proposition}
\label{71}
$\DAHA$ can equivalently be described as the algebra generated by
$T_1$, $Y$, $Y^{-1}$, $Z$, $Z^{-1}$ with relations
$YY^{-1}=1=Y^{-1}Y$, $ZZ^{-1}=1=Z^{-1}Z$ and
\begin{gather}
T_1^2=-(ab+1)T_1-ab,\nonumber\\
T_1Z = Z^{-1}T_1+(ab+1)Z^{-1}-(a+b),\nonumber\\
T_1Z^{-1}= ZT_1-(ab+1)Z^{-1}+(a+b),\nonumber\\
T_1Y= q^{-1}abcd Y^{-1}T_1-(ab+1)Y+ab(1+q^{-1}cd),\nonumber\\
T_1Y^{-1}= q(abcd)^{-1}YT_1+q(abcd)^{-1}(1+ab)Y-q(cd)^{-1}(1+q^{-1}cd),\nonumber\\
YZ= qZY+(1+ab)cd\,Z^{-1}Y^{-1}T_1
-(a+b)cd\,Y^{-1}T_1
-(1+q^{-1}cd)Z^{-1}T_1\nonumber\\
\phantom{YZ=}{}
-(1-q)(1+ab)(1+q^{-1}cd)Z^{-1}
+(c+d)T_1
+(1-q)(a+b)(1+q^{-1}cd),\nonumber\\
YZ^{-1}= q^{-1}Z^{-1}Y
-q^{-2}(1+ab)cd\,Z^{-1}Y^{-1}T_1
+q^{-2}(a+b)cd\,Y^{-1}T_1\nonumber\\
\phantom{YZ^{-1}=}{}
+q^{-1}(1+q^{-1}cd)Z^{-1}T_1
-q^{-1}(c+d)T_1,\nonumber\\
Y^{-1}Z=q^{-1}ZY^{-1}-q(ab)^{-1}(1+ab)Z^{-1}Y^{-1}T_1
+(ab)^{-1}(a+b)Y^{-1}T_1\nonumber\\
\phantom{Y^{-1}Z=}{}
+q(abcd)^{-1}(1+q^{-1}cd)Z^{-1}T_1
+q(abcd)^{-1}(1-q)(1+ab)(1+q^{-1}cd)Z^{-1}\nonumber\\
\phantom{Y^{-1}Z=}{}
-(abcd)^{-1}(c+d)T_1
-(abcd)^{-1}(1-q)(1+ab)(c+d),\nonumber\\
Y^{-1}Z^{-1}= qZ^{-1}Y^{-1}+q(ab)^{-1}(1+ab)Z^{-1}Y^{-1}T_1
-(ab)^{-1}(a+b)Y^{-1}T_1\nonumber\\
\phantom{Y^{-1}Z^{-1}=}{}
-q^2(abcd)^{-1}(1+q^{-1}cd)Z^{-1}T_1
+q(abcd)^{-1}(c+d)T_1.
\label{70}
\end{gather}
$\DAHA$ is spanned by the elements $Z^mY^nT_1^i$, where
$m,n\in\ZZ$, $i=0,1$.
\end{Proposition}
\begin{proof}
First we start with relations \eqref{64}--\eqref{69}.
Then \eqref{64}, \eqref{65} give \eqref{38}, \eqref{39}. Next put
$Y:=T_1T_0$, $Y^{-1}:=T_0^{-1}T_1^{-1}$.
Then verify relations \eqref{70} from relations \eqref{64}--\eqref{69},
most conveniently with the aid of computer algebra package, for instance
by using \cite{13}.

Conversely we start with relations \eqref{70}. Then the f\/irst of these
relations gives \eqref{38}. Put $T_0:=T_1^{-1}Y$.
Then verify relations \eqref{64}--\eqref{69} from relations \eqref{70},
where again computer algebra may be used.

The last statement follows from the PBW-type structure of the
relations \eqref{70}. Observe that by the f\/irst f\/ive relations together
with the trivial relations, every word in $T_1$, $Y$, $Y^{-1}$, $Z$, $Z^{-1}$ can
be written as a linear combination of words with at most one occurrence
of $T_1$ in each word and only on the right, and with no substrings
$YY^{-1}$, $Y^{-1}Y$, $ZZ^{-1}$, $Z^{-1}Z$, and with no more occurrences of
$Y$, $Y^{-1}$, $Z$, $Z^{-1}$ in each word than in the original word.
If in one of these terms there are misplacements ($Y$ or $Y^{-1}$ before
$Z$ or $Z^{-1}$) then apply one of the last four relations followed by the
previous step in order to reduce the number of misplacements.
\end{proof}

\begin{Theorem}
The basic representation \eqref{14}--\eqref{16} of $\DAHA$ is faithful.
A basis of $\DAHA$ is provided by the elements $Z^mY^nT_1^i$, where
$m,n\in\ZZ$, $i=0,1$.
\end{Theorem}
\begin{proof}
For $j>0$ we have
\begin{gather}
Z^mY^n\,E_{-j}=q^{-jn}\,z^{m-j}+\cdots+\const z^{m+j-1},\nonumber\\
Z^mY^nT_1\,E_{-j}=\const z^{m-j}+\cdots-ab(q^{j-1}abcd)^n\,z^{m+j},\label{72}\\
Z^mY^nT_1^{-1}\,E_{-j}=\const z^{m-j}+\cdots+(q^{j-1}abcd)^n\,z^{m+j}.\nonumber
\end{gather}
This follows from \eqref{25}--\eqref{28}, \eqref{52} and \eqref{38}.
Suppose that some linear combination
\begin{gather}
\sum_{m,n}a_{m,n} Z^mY^n+\sum_{m,n}b_{m,n}Z^mY^nT_1
\label{73}
\end{gather}
acts as the zero operator in the basic representation, while not all
coef\/f\/icients $a_{m,n}$, $b_{m,n}$ are zero.
Then there is a maximal $r$ for which $a_{r,n}$ or $b_{r,n}$ is nonzero
for some~$n$. If $b_{r,n}\ne0$ for some~$n$ then let the operator
\eqref{73} act on $E_{-j}$.
By \eqref{72} we have
that for all $j\ge1$
\begin{gather*}
\sum_n b_{r,n}(q^{j-1}abcd)^n\,z^{r+j}=0,\qquad
{\rm hence}\qquad \sum_n b_{r,n}(q^{j-1}abcd)^n=0.
\end{gather*}
By assumption \eqref{21} we see that $\sum_n b_{r,n} w^n=0$.
Hence $b_{r,n}=0$ for all $n$, which is a contradiction.

So $a_{r,n}\ne0$ for some $n$. Let the operator
\eqref{73} act on $T_1^{-1}E_{-j}$.
By \eqref{72} we have
that for all $j\ge1$
\begin{gather*}
\sum_n a_{r,n}(q^{j-1}abcd)^n\,z^{r+j}=0,\qquad
{\rm hence}\qquad \sum_n a_{r,n}(q^{j-1}abcd)^n=0.
\end{gather*}
Again we arrive at the contradiction that $a_{r,n}=0$ for all $n$.
\end{proof}

\section[The embedding of a central extension of $AW(3,Q_0)$ in $\tilde{\mathfrak{H}}$]{The embedding of a central extension of $\boldsymbol{AW(3,Q_0)}$ in $\boldsymbol{\DAHA}$}

Let us now examine whether the representation \eqref{4} of $AW(3)$ on
$\Asym$ extends to a representation on $\FSA$ if we let $K_0$ act as
$D$ instead of $\Dsym$.
It will turn out that this is only true
for certain specializations of $a$, $b$, $c$, $d$, but that a suitable central
extension $\AW$ of $AW(3)$ involving $T_1$
will realize what we desire.

\begin{Definition}\rm
$\AW$ is the algebra generated by $K_0$, $K_1$, $T_1$ with relations
\begin{gather}
\label{43}
T_1K_0=K_0T_1,\qquad
T_1K_1=K_1T_1,\qquad
(T_1+ab)(T_1+1)=0,\\
(q+q^{-1})K_1K_0K_1-K_1^2K_0-K_0K_1^2\nonumber\\
\qquad{}{}=B\,K_1+ C_0\,K_0+D_0+E\,K_1(T_1+ab)
+F_0(T_1+ab),\label{44}\\
(q+q^{-1})K_0K_1K_0-K_0^2K_1-K_1K_0^2\nonumber\\
\qquad{}=B\,K_0+C_1\,K_1+D_1+E\,K_0(T_1+ab)
+F_1(T_1+ab),\label{45}
\end{gather}
where the structure constants are given by \eqref{42} together with
\begin{gather}
E:=-q^{-2}(1-q)^3(c+d),\nonumber\\
F_0:=q^{-3}(1-q)^3(1+q)(cd+q),\label{46}\\
F_1:=q^{-3}(1-q)^3(1+q)(a+b)cd.\nonumber
\end{gather}
\end{Definition}

It can be shown that the following adaptation
of \eqref{80} is a Casimir operator for $\AW$, commuting with $K_0$, $K_1$, $T_1$:
\begin{gather}
\Q:=(K_1K_0)^2-(q^2+1+q^{-2})K_0(K_1K_0)K_1+(q+q^{-1})K_0^2K_1^2
\nonumber\\
\phantom{\Q:=}{}+(q+q^{-1})(C_0K_0^2+C_1K_1^2)+\bigl(B+E(T_1+ab)\bigr)\bigl((q+1+q^{-1})K_0K_1+K_1K_0\bigr)
\nonumber\\
\phantom{\Q:=}{}+(q+1+q^{-1})\bigl(D_0+F_0(T_1+ab)\bigr)K_0+(q+1+q^{-1})\bigl(D_1+F_1(T_1+ab)\bigr)K_1\nonumber\\
\phantom{\Q:=}{}+G(T_1+ab),
\label{87}
\end{gather}
where
\begin{gather}
G:=-q^{-4}(1-q)^3\Bigl((a+b)(c+d)\bigl(cd(q^2+1)+q\bigr)
-q(ab+1)\bigl((c^2+d^2)(q+1)-cd\bigr)\nonumber\\
\phantom{G:=}{}+(cd+e_4)(q^2+1)+(e_2+e_4-ab)q^3\Bigr).
\label{90}
\end{gather}
Let $\AWQ$ be the algebra generated by $K_0$, $K_1$, $T_1$ with relations
\eqref{43}--\eqref{45} and additional relation
\begin{gather}
\Q=Q_0,
\label{91}
\end{gather}
where $\Q$ is given by \eqref{87} and $Q_0$ by \eqref{81}.

\begin{Theorem}
There is a representation of the algebra $\AWQ$ on the space $\FSA$
of Laurent polynomials $f[z]$ such that $K_0$ acts as $D$, $K_1$ acts
by multiplication by $z+z^{-1}$, and the action of $T_1$ is given by
\eqref{15}.
This representation is faithful.
\end{Theorem}
\begin{proof}
It follows by straightforward computation, possibly using computer algebra,
that this is a representation of $\AWQ$.
In the same way as for Lemma \ref{55} it can be shown that $\AWQ$ is spanned
by the elements
\begin{gather}
K_0^n(K_1K_0)^iK_1^mT_1^j\qquad(m,n=0,1,2,\ldots,\ \ i,j=0,1).
\label{74}
\end{gather}
Now we will prove that the representation is faithful.
Suppose that for certain coef\/f\/icients $a_{k,l}$, $b_{k,l}$, $c_{k,l}$, $d_{k,l}$
we have
\begin{gather}
\sum_{k,l}a_{k,l}\,D^l\,(Z+Z^{-1})^k+
\sum_{k,l}b_{k,l}\,D^{l-1}\,(Z+Z^{-1})\,D\,(Z+Z^{-1})^{k-1}\nonumber\\
+\Bigg(\sum_{k,l}c_{k,l}\,D^l\,(Z+Z^{-1})^k+
\sum_{k,l}d_{k,l}\,D^{l-1}\,(Z+Z^{-1})\,D\,(Z+Z^{-1})^{k-1}\Bigg)(T_1+ab)=0
\label{75}
\end{gather}
while acting on $\FSA$.
Then, since $T_1P_j=-ab P_j$ (see \eqref{51}),
we have for all $j\ge0$ that
\begin{gather*}
\sum_{k,l}a_{k,l}\,\Dsym^l\,(Z+Z^{-1})^k\,P_j[z]+
\sum_{k,l}b_{k,l}\,\Dsym^{l-1}\,(Z+Z^{-1})\,D\,(Z+Z^{-1})^{k-1}\,P_j[z]=0.
\end{gather*}
Then by the proof of Theorem \ref{76} if follows that all coef\/f\/icients
$a_{k,l}$, $b_{k,l}$ vanish.

It follows from \eqref{49} and \eqref{51} that
$(T_1+ab)E_{-n}=-ab Q_n$ (also if $ab=1$). Hence,
if we let~\eqref{75}, with vanishing $a_{k,l}$, $b_{k,l}$, act on $E_{-j}[z]$,
and divide by $-ab$, then:
\begin{gather*}
\Bigg(\sum_{k,l}c_{k,l}\,D^l\,(Z+Z^{-1})^k+
\sum_{k,l}d_{k,l}\,D^{l-1}\,(Z+Z^{-1})\,D\,(Z+Z^{-1})^{k-1}\Bigg)Q_j[z]=0.
\end{gather*}
{}From \eqref{31} we see that the three-term recurrence relation
\eqref{59} for $P_n[z]$
has an analogue for~$Q_n[z]$:
\begin{gather*}
(z+z^{-1})Q_n[z]=Q_{n+1}[z]+\tilde\be_n Q_n[z]
+\tilde\ga_n Q_{n-1}[z]\qquad(n\ge2),
\end{gather*}
where $\tilde\be_n$ and $\tilde\ga_n$ are obtained from the
corresponding $\be_n$ and $\ga_n$ (\eqref{60} and \eqref{61})
by replacing $a$, $b$, $n$ by $qa$, $qb$, $n-1$, respectively.
Hence \eqref{84} remains valid if we replace each $P$ by $Q$.
Again, similarly as in the proof of Theorem \ref{76},
it follows that all coef\/f\/icients
$c_{k,l}$, $d_{k,l}$ vanish.
\end{proof}

\begin{Corollary}
\label{77}
The algebra $\AWQ$ can be isomorphically embedded into $\DAHA$ by the mapping
\begin{gather}
K_0\mapsto Y+q^{-1}abcd Y^{-1},\qquad
K_1\mapsto Z+Z^{-1},\qquad
T_1\mapsto T_1.
\label{78}
\end{gather}
\end{Corollary}
\begin{proof}
The embedding is valid for $\AWQ$ and $\DAHA$ acting on $\FSA$.
Now use the faithfulness of the representations of $\AWQ$ and $\DAHA$ on $\FSA$.
\end{proof}
\begin{Remark}\rm
By Corollary \ref{77} the relations \eqref{43}--\eqref{45} and \eqref{91}
are valid
identities in $\DAHA$ after substitution by \eqref{78}. These identities
can also be immediately verif\/ied within $\DAHA$, for instance by usage
of the package \cite{13}.
\end{Remark}
\begin{Remark}\rm
If $a$, $b$, $c$, $d$ are such that $E,F_0,F_1=0$ in \eqref{46} then we have already
a homomorphism of the original algebra $AW(3)$ into $\DAHA$ under the
substitutions $K_0:=D$, $K_1:=Z+Z^{-1}$
in \eqref{1}, \eqref{2}.
This is the case
if\/f $c=-d=q^\half$ (or $-q^\half$) and $a=-b$. For these parameters the
Askey--Wilson polynomials become the continuous $q$-ultraspherical
polynomials
(see \cite[(7.5.25), (7.5.34)]{2}):
\begin{gather*}
P_n\big[z;a,-a,q^\half,-q^\half\mid q\big]=\const
C_n\big(\thalf(z+z^{-1});a^2\mid q^2\big).
\end{gather*}
However, for these specializations of
$a$, $b$, $c$, $d$ we see from \eqref{80} and \eqref{87}
that $\Q$ still slightly dif\/fers from $Q$:
it is obtained from $Q$ by adding the term
$(q^{-1}-q)^3(1-a^2)(T_1-a^2)$. So $\AWQ$ then still
dif\/fers from $AW(3,Q_0)$.

For such $a$, $b$, $c$, $d$ the operator $T_0$ acting on $\FSA$
(formula \eqref{16}) simplif\/ies to
$(T_0f)[z]=f[qz^{-1}]$. We then have the specialization of parameters in
$\DAHA$ to the one-parameter double af\/f\/ine Hecke algebra of type $A_1$
(see \cite[\S~6.1--6.3]{5}). Explicit formulas for the non-symmetric
$q$-ultraspherical polynomials become much nicer than in the general
four-parameter Askey--Wilson case, see \cite[(6.2.7), (6.2.8)]{5}.
\end{Remark}
\subsection*{Acknowledgements}
I am very much indebted to an anonymous referee who pointed out errors
in the proof of an earlier version of Theorem \ref{76}. I also thank him
for suggesting a further simplification in my proof of the corrected theorem.
I thank
Siddhartha Sahi for making available to me a draft of his paper \cite{12}
in an early stage. I also thank Jasper Stokman for helpful comments.

\pdfbookmark[1]{References}{ref}
\LastPageEnding
\end{document}